\documentclass[12pt,reqno]{article}
\usepackage{}
\usepackage{amsfonts,amsmath,amssymb,latexsym}
\usepackage{mathrsfs}
\usepackage{CJK}
\usepackage{fancyhdr}
\usepackage{hyperref}
\usepackage{anysize}
\usepackage{graphicx,epsfig} 
\usepackage{amsthm}
\usepackage[perpage,symbol*]{footmisc}
\usepackage{cite}
\usepackage{xcolor}

\newtheorem{thm}{Theorem}[section]
\newtheorem{defi}[thm]{Definition}
\newtheorem{prop}[thm]{Proposition}
\newtheorem{lem}[thm]{Lemma}

\newtheorem{rem}[thm]{Remark}

\makeatletter \@addtoreset{equation}{section} \makeatother
\def\beq{\begin{equation}}
\def\eeq{\end{equation}}
\def\endproof{$\hfill\Box$}
\begin{document}
\baselineskip=20pt  \hoffset=-3cm \voffset=0cm \oddsidemargin=3.2cm
\evensidemargin=3.2cm \thispagestyle{empty}\vspace{10cm}
\title{\textbf{Homoclinic orbits for a class of nonperiodic
Hamiltonian systems with some twisted conditions }}
\author{\Large Qi Wang $^{{\rm a},\,*}$ ,$\quad$  Qingye Zhang
$^{{\rm b}}$}
\date{} \maketitle

\begin{center}
\it\scriptsize ${}^{\rm a}$Institute of Contemporary Mathematics, School of Mathematics and Information Science,\\ Henan
University, Kaifeng 475000, PR China\\
${}^{\rm b}$Department of Mathematics,
Jiangxi Normal University, Nanchang 330022, PR China
\end{center}

\footnotetext[0]{$^*$Corresponding author. Partially supported by NNSF (10901118,11126154).}

\footnotetext[0]{\footnotesize\;{\it E-mail address}:
mathwq@henu.edu.cn. (Q. Wang), qingyezhang@gmail.com (Q. Zhang).}

\noindent
{\bf Abstract:} {\small In this paper, by the Masolv index theory, we will study the existence and multiplicity of homoclinic orbits for a class of asymptotically linear nonperiodic Hamiltonian systems with some twisted conditions on the Hamiltonian functions.}


\section{Introduction and main results}

Consider the following
first order non-autonomous Hamiltonian systems
\[
\dot{z}=JH_z(t,z),\eqno{\rm (HS)}
\]
where $z:\mathbb{R}\to \mathbb{R}^{2N}$, $J=\left(\begin{matrix} 0 &-I_N\\I_N &0\end{matrix}\right)$,
 $H\in C^1(\mathbb{R}\times\mathbb{R}^{2N},\mathbb{R})$ and $\nabla_zH(t,z)$ denotes the gradient of $H(t,z)$ with respect to $z$. As usual we say that a nonzero solution
$z(t)$ of (HS) is homoclinic (to 0) if $z(t)\to 0$  as $|t|\to\infty$.

As a special case of dynamical systems, Hamiltonian systems are very important in the study of gas dynamics, fluid
mechanics, relativistic mechanics and nuclear physics. While it is well known that homoclinic solutions play an important role in analyzing the chaos of Hamiltonian systems. If a system has the transversely intersected homoclinic solutions, then it must be chaotic. If it has the smoothly connected homoclinic solutions,
then it cannot stand the perturbation, its perturbed system probably produces chaotic phenomena. Therefore, it is of practical
importance and mathematical significance to consider the existence of homoclinic solutions of Hamiltonian systems emanating from 0.

In the last years, the existence and multiplicity of homoclinic orbits for the first order system
(HS) were studied extensively by means of critical point theory, and many results were obtained
under the assumption that $H(t,z)$ depends periodically on $t$ (see, e.g., \cite{AS,CES,D,DG,DW,HW,S92,S93,SZ,T}). Without assumptions of periodicity the problem is quite different in nature and there is not much work done so far. To the best of our knowledge, the authors in \cite{Ding-1995-JMAA} firstly obtained the existence of homoclinic orbits for a class of first order systems without any periodicity on the Hamiltonian function. After this, there were a few papers dealing with the the existence and multiplicity of homoclinic orbits for the first order system (HS) in this situation (see, e.g., \cite{DJ-JDE-2007,DL,SCN}).

In the present paper, with the Maslov index theory of homoclinic orbits introduced by Chen and Hu in \cite{Chen-Hu-2007}, we will study the existence and multiplicity of homoclinic orbits for (HS) without any periodicity on the Hamiltonian function. To the best of the author's knowledge, the Maslov index theory of homoclinic orbits is the first time to be used to study the existence of homoclinic solutions. We are mainly interested in the Hamiltonian functions of the form
\begin{equation}\label{Hamiltonian Function}
H(t,z)=-L(t)z\cdot z+R(t,z),
\end{equation}
where $L$ is an $2N\times 2N$ symmetric matrix valued function.
We assume that\\
($L_1$)$L\in C(\mathbb{R},\mathbb{R}^{N^2})$, and there are $\alpha,\;c>0$, $t_0\geq 0$ and a constant matrix $P$, satisfying
$$
PL(t)-c|t|^{\alpha}I_{2N}\geq 0,\;\forall |t|\geq t_0,
$$
where $I_{2N}$ is the identity map on $\mathbb{R}^{2N}$ and  for a $2N\times 2N$ matrix $M$, we say $M\geq 0$ if and only if
$$
\displaystyle\inf_{\xi\in\mathbb{R}^{2N},|\xi|=1}M\xi\cdot\xi\geq 0.
$$
In ($L_1$), if $P=\left(\begin{matrix} 0 &I_N\\I_N &0\end{matrix}\right)$, then ($L_1$) is similar to  the condition ($R_0$) in \cite{DJ-JDE-2007}. But the restrictions on $R(t,z)$ will be  different from \cite{DJ-JDE-2007}, and we will give some examples in Remark \ref{examples}.
If $P=\pm I_{2N}$ or $P=\left(\begin{matrix} I_{N+m} &0\\0 &-I_{N-m}\end{matrix}\right)$ in condition ($L_1$), for examples, it's quite different from the existing  results as authors known. In short, condition ($L_1$) means that the eigenvalues of $L(t)$ will tend to $\pm\infty$ with the speed no less than $|t|^{\alpha}$. But ($L_1$) does not contain all of these cases. For examples, let $N=1$ and $L(t)=|t|^\alpha\left(\begin{matrix} \cos 2t &\sin 2t\\\sin 2t &-\cos 2t\end{matrix}\right)$, we have the eigenvalues of $L(t)$ are $\pm |t|^\alpha$, but there is no constant matrix $P$ satisfying $PL(t)-c|t|^{\alpha}I_{2N}\geq 0,\;\forall |t|\geq t_0$.

Denote by $\tilde{F}$ the self-adjoint operator $-\mathcal {J}\frac{d}{dt}+L(t)$ on $ L^2\equiv L^2(\mathbb{R},\mathbb{R}^{2N})$, with domain
$D(\tilde{F})=H^1(\mathbb{R},\mathbb{R}^{2N})$ if $L(t)$ is bounded and $D(\tilde{F})\subset H^1(\mathbb{R},\mathbb{R}^{2N})$ if $L(t)$ is unbounded. Let  $|\tilde{F}|$  be  the absolute value of $\tilde{F}$, and
 $|\tilde{F}|^{1/2}$ be the square root of $|\tilde{F}|$. $D(\tilde{F})$ is a Hilbert space equipped with the norm
\beq
||z||_{\tilde{F}}=||(I+|\tilde{F}|)z||_{L^2},\;\forall z\in D(\tilde{F}).
\eeq
 Let
$E=D(|\tilde{F}|^{1/2})$, and define
on $E$ the inner product and norm by
$$(u,v)_E=(|\tilde{F}|^{1/2}u,|\tilde{F}|^{1/2}v)_2+(u,v)_2,$$
$$\|u\|_E=(u,u)^{1/2}_E,$$
where $(\cdot,\cdot)_{L^2}$ denotes the usual inner product on $L^2(\mathbb{R},\mathbb{R}^{2N})$.
Then $E$ is a Hilbert space. It is easy to see that $E$ is
continuously embedded in $H^{1/2}(\mathbb{R}, \mathbb{R}^{2N})$,
and we further have the following lemma.

\begin{lem}\label{Lp-embedding}
 Suppose that $L$ satisfies $(L_1)$. Then  $E$ is compactly embedded  in
$L^p(\mathbb{R},\mathbb{R}^{2N})$ with the usual norm $\|\cdot\|_{L^p}$ for any $1\leq p\in (\frac{2}{1+\alpha},\infty)$.
\end{lem}

This lemma is similar to Lemma 2.1-2.3 in \cite{Ding-1995-JMAA}, and  we will prove it in Section 3.
Define the quadratic form $\mathcal {Q}$ on $E$ by
\begin{equation}
\mathcal {Q}(u,v)=\int_\mathbb{R}((-J\dot u, v)+(L(t)u,v))dt,\;\,\forall\,  u,v\in E.
\end{equation}
It's easy to check that $\mathcal {Q}(u,v)$ is a bounded quadratic form on $E$ and hence there exists a unique bounded self-adjoint operator $F:E\to E$ such that
\begin{equation}\label{operator-F}
 (Fu,v)_E=\mathcal {Q}(u,v),\,\forall\, u,v\in E.
\end{equation}
Besides, define a linear operator $K:L^2\to E$ by
 \begin{equation}\label{operator-K}
 (Ku,v)_E=(u,v)_{L^2},\,\forall\,  u\in L^2,v\in E.
\end{equation}
In view of Lemma \ref{Lp-embedding}, we know that $F$ is a Fredholm operator and $K$ is a compact operator.

Denote by   $\mathcal{B}$ the set of all uniformly  bounded  symmetric $2N\times 2N$ matric functions.
That is to say $B\in\mathcal{B}$ if and only if $B^T(t)=B(t)$ for all $t\in\mathbb{R}$ and
 $B(t)$ is uniformly bounded  in $t$ as the operator on $\mathbb{R}^{2N}$.
For any $B\in \mathcal{B}$, it is easy to see $B$ determines a  bounded self-adjoint operator  on $L^2$, by $z(t)\mapsto B(t)z(t)$, for any
$z\in L^2$, we still denote this operator by $B$,
then $KB:E\subset L^2\to E$ is a self-adjoint compact operator on $E$ and  satisfies
 \begin{equation}
 (KB u,v)_E=(Bu,v)_{L^2}, \,\forall\,  u,v\in E.
\end{equation}

 Before presenting the conditions on $R(t,z)$, we need the concept of Maslov index for homoclinic orbits introduced by Chen and Hu in \cite{Chen-Hu-2007} which  is equivalent to the  relative Morse index.
 We will give a brief introduction of it by Definition \ref{relative fredholm index}, where for any $B\in\mathcal{B}$, we denote  the associated index pair by ($\mu_F(KB)$, $\upsilon_F(KB)$).

Now we can present the conditions on $R(t,z)$ as follows. For notational simplicity, we set $B_0(t)=\nabla^2_zR(t,0)$, and in what follows
the letter $c$ will be repeatedly used to denote various positive constants whose exact value is irrelevant.
Besides, for two $2N\times 2N$ symmetric matrices $M_1$ and $M_2$, $M_1\le M_2$ means that $M_2-M_1$ is semi-positive definite.

\begin{list}{}
{\setlength{\topsep}{1ex} \setlength{\itemsep}{0ex}
\setlength{\leftmargin}{2.6em}
 }
 \item[$(R_1)$] $R\in C^2(\mathbb{R}\times\mathbb{R}^{2N},\mathbb{R})$, and there exists a constant $c>0$ such that
\[
|\nabla^2_zR(t,z)|\leq c,\;\,\forall\,  (t,z)\in \mathbb{R}\times\mathbb{R}^{2N}.
\]

\item[$(R_0)$] $\nabla_zR(t,0)\equiv 0$ and  $B_0\in \mathcal{B}$.

\item[$(R_\infty)$] There exists some $R_0>0$ and  continuous symmetric matrix functions $B_1,B_2\in\mathcal{B}$ with $\mu_F(KB_1)=\mu_F(KB_2)$ and $\upsilon_F(KB_2)=0$ such that
\[
B_1(t)\leq\nabla^2_zR(t,z)\leq  B_2(t),\;\,\forall\,  t\in\mathbb{R},\; |z|>R_0.
\]
\end{list}
 Then we have our first result.
\begin{thm}\label{main-result3}
Assume $(L_1)$, $(R_1)$,  $(R_0)$ and $(R_\infty)$ hold.
If
\[\mu_F(KB_1)\not\in [\mu_F(KB_0), \mu_F(KB_0)+\upsilon_F(KB_0)],\]
 then {\rm (HS)} has at least one nontrivial homoclinic orbit. Moreover, if $\upsilon_F(KB_0)=0$ and $|\mu_F(KB_1)-\mu_F(KB_0)|\geq N$, the problem possesses
 at least two nontrivial homoclinic orbits.
\end{thm}
Condition $(R_\infty)$ is a two side pinching condition near the infinity,  we can relax $(R_\infty)$ to condition $(R^{\pm}_\infty)$ as follows.
\begin{list}{}
{\setlength{\topsep}{1ex} \setlength{\itemsep}{0ex}
\setlength{\leftmargin}{2.6em}
 }

\item[$(R^{\pm}_\infty)$] There exists some $R_0>0$ and a continuous symmetric matrix function $B_\infty\in\mathcal{B}$
with $\upsilon_F(KB_\infty)=0$ such that
\[
\pm \nabla^2_zR(t,z)\geq \pm B_\infty(t),\;\,\forall\,  t\in\mathbb{R},\; |z|>R_0.
\]
\end{list}
Then we have the following results.
\begin{thm}\label{main-result1}Assume $(L_1)$, $(R_1)$, $(R_0)$, $(R^{+}_\infty)(or\, (R^{-}_\infty))$ and  $\upsilon_F(KB_0)=0$  hold.
If $\mu_F(KB_\infty)\geq \mu_F(KB_0)+2\;(or \;\, \mu_F(KB_\infty)\leq \mu_F(KB_0)-2)$, then {\rm (HS)} has at least
one nontrivial homoclinic orbit.
\end{thm}
\begin{thm}\label{main-result2}
Suppose that $(L_1)$, $(R_1)$, $(R_0)$,  $(R^{+}_\infty)\;(or\, (R^{-}_\infty))$ and  $\upsilon_F(KB_0)=0$  are satisfied. If in addition, $R$ is even in $z$ and $\mu_F(KB_\infty)\geq \mu_F(KB_0)+2 \;(or\; \, \mu_F(KB_\infty)\leq \mu_F(KB_0)-2)$, then {\rm (HS)} has at least
$|\mu_F(KB_\infty)-\mu_F(KB_0)|-1$ pairs of nontrivial homoclinic orbits.
\end{thm}
\begin{rem}\label{examples}
{\rm  Lemma \ref{Lp-embedding} shows that $\sigma(A)$, the spectrum of $A$, consists of eigenvalues numbered by (counted in their multiplicities):
\[
\cdots\leq \lambda_{-2}\leq\lambda_{-1}\leq 0<\lambda_1\leq\lambda_2\leq\cdots
\]
with $\lambda_{\pm k}\to\pm \infty$ as $k\to \infty$. Let $B_0(t)\equiv B_0$ and $B_\infty(t)\equiv B_\infty$, with the constants $B_0, B_\infty$
 satisfying $\lambda_l<B_0<\lambda_{l+1}$,
and $\lambda_{l+i}<B_\infty<\lambda_{l+i+1}$ for some $l\in\mathbb{Z}$ and $i\geq 1$ (or $i\leq -1$). Define
\[
R(t,z)=\delta(|z|)\frac{1}{2}B_0|z|^2+(1-\delta(|z|))\frac{1}{2}B_\infty|z|^2,
\]
where $\delta$ is a smooth cutoff function satisfying $\delta(|z|)=\left\{\begin{array}{ll}1,\;|z|< 1,\\0,\; |z|>2. \end{array}\right.$ By Proposition \ref{index-sp} below, it is easy to verify $R$ satisfies all the conditions in Theorem \ref{main-result3}. Furthermore,  let  the constant $B_\infty$ satisfying
 $\lambda_{l+i}<B_\infty<\lambda_{l+i+1}$ for some $l\in\mathbb{Z}$ and $i\geq 2$ (or $i\leq-2$). Define
\[
R(t,z)=\delta(|z|)\frac{1}{2}B_0|z|^2+(1-\delta(|z|))\frac{1}{2}B_\infty|z|^2,
\]
Then $R$ satisfies all the conditions in Theorem \ref{main-result1} and Theorem \ref{main-result2}. However, it is easy to see that some conditions of the main results in \cite{Ding-1995-JMAA,DJ-JDE-2007,DL,SCN} does not hold for these examples.
}
\end{rem}

\begin{rem}\label{rem-1.4}
{ \rm
 Note that the assumption $\upsilon_F(KB_\infty)=0$ in $(R^{\pm}_\infty)$ is not essential  for our main results. For the case of  $(R^{+}_\infty)$ with $\upsilon_F(KB_\infty)\neq 0$,  let $\widetilde{B}_\infty=B_\infty-\varepsilon I_{2N}$ with $\varepsilon >0$ small enough, where $I_{2N}$ is the identity map on $\mathbb{R}^{2N}$, then  $\mu_F(K\widetilde{B}_\infty)=\mu_F(KB_\infty)$ and $\upsilon_F(K\widetilde{B}_\infty)= 0$, and hence $(R^{+}_\infty)$ holds for $\widetilde{B}_\infty$. Therefore Theorems \ref{main-result1} and \ref{main-result2} still hold in this case. While for the case of $(R^{-}_\infty)$ with $\upsilon_F(KB_\infty)\neq 0$, if we replace $\mu_F(KB_\infty)$ by
$ \mu_F(KB_\infty)+\upsilon_F(KB_\infty)$ in Theorems \ref{main-result1} and \ref{main-result2}, then similar results hold. Indeed, let $\widetilde{B}_\infty=B_\infty+\varepsilon I_{2N}$ with $\varepsilon>0$ small enough such that $\mu_F(K\widetilde{B}_\infty)=\mu_F(KB_\infty)+\upsilon_F(KB_\infty)$ and $\upsilon_F(K\widetilde{B}_\infty)= 0$, then this case is also reduced to the case of $(W^{-}_\infty)$ for $\widetilde{B}_\infty$ with $\upsilon_F(K\widetilde{B}_\infty)= 0$.}
\end{rem}

\section{Preliminaries}

In this section, we recall  the definition of relative Morse index, saddle point reduction, and  give the relationship between them. For this propose, the notion of spectral flow will be used.
 \subsection{Relative Morse index}
Let $\mathcal{H}$ be a separable Hilbert space, for any self-adjoint operator $A$ on $\mathcal{H}$, there is a unique $A$-invariant orthogonal splitting
\begin{equation}\label{invarsplit}
\mathcal{H}=\mathcal{H}^+(A)\oplus \mathcal{H}^-(A)\oplus \mathcal{H}^0(A),
\end{equation}
where $\mathcal{H}^0(A)$ is the null space of $A$, $A$ is positive definite on $\mathcal{H}^+(A)$ and negative definite on $\mathcal{H}^-(A)$, and $P_A$ denotes the orthogonal projection from $\mathcal{H}$ to $\mathcal{H}^-(A)$. For any bounded self-adjoint Fredholm operator $\mathcal{F}$ and a compact self-adjoint operator $\mathcal{T}$ on $\mathcal{H}$,
$P_\mathcal{F}-P_{\mathcal{F}-\mathcal{T}}$ is compact (see  Lemma 2.7 of \cite{Zhu-Long-cam}), where $P_\mathcal{F}:\mathcal{H}\to \mathcal{H}^-(\mathcal{F})$ and $P_{\mathcal{F}-\mathcal{T}}:\mathcal{H}\to \mathcal{H}^-(\mathcal{F}-\mathcal{T})$ are the respective projections.
Then by Fredholm operator theory, $P_\mathcal{F}|_{\mathcal{H}^-{(\mathcal{F}-\mathcal{T})}}:\mathcal{H}^-{(\mathcal{F}-\mathcal{T})}\to \mathcal{H}^-{(\mathcal{F})}$ is a Fredholm operator. Here and in the sequel, we denote by $\rm ind (\cdot)$ the Fredholm index of a Fredholm operator.
\begin{defi}\label{relative fredholm index}
For any bounded self-adjoint Fredholm operator $\mathcal{F}$ and a compact self-adjoint operator $\mathcal{T}$ on $\mathcal{H}$, the relative Morse index pair $(\mu_\mathcal{F}(\mathcal{T}), \upsilon_\mathcal{F}(\mathcal{T}))$ is defined by
\begin{equation}
\mu_\mathcal{F}(\mathcal{T})=\rm ind(P_\mathcal{F}|_{\mathcal{H}^-{(\mathcal{F}-\mathcal{T})}}).
\end{equation}
and
 \begin{equation}
 \upsilon_\mathcal{F}(\mathcal{T})=\dim \mathcal{H}^0(\mathcal{F}-\mathcal{T}).
 \end{equation}
  \end{defi}
\subsection{Saddle point reduction}\label{s-p-r}
In this subsection, we describe the saddle point reduction  in \cite{Amann,Chang-1993,Long3}.
Recall that $\mathcal{H}$ is a real Hilbert space, and  A is a self-adjoint operator with domain $D(A)\subset \mathcal{H}$. Let $\Phi\in C^1(\mathcal{H},\mathbb{R})$, with $\Phi'(\theta)=0$. Assume that
\begin{list}{}
{\setlength{\topsep}{1ex} \setlength{\itemsep}{0ex}
\setlength{\leftmargin}{2.6em}
 }
\item[$(1)$]
 There exist real numbers $\alpha<\beta$ such that $\alpha,\beta\notin\sigma(A)$, and that $\sigma(A)\cap[\alpha,\beta]$ consists of at most finitely many eigenvalues of finite multiplicities.

\item[$(2)$]
$\Phi'$ is Gateaux differentiable in $\mathcal{H}$, which satisfies
\[
\|d\Phi'(u)-\frac{\alpha+\beta}{2}I\|\leq\frac{\beta-\alpha}{2},\;\forall u\in \mathcal{H}.
\]
Without loss of generality, we may assume $\alpha=-\beta, \beta>0$.

\item[$(3)$]
 $\Phi\in C^2(V,\mathbb{R})$,  $V=D(|A|^{1/2})$, with the norm
\[
\|z\|_V=(\||A|^{1/2}z\|^2_\mathcal{H}+\varepsilon^2\|z\|^2_\mathcal{H})^{1/2},
\]
where $\varepsilon>0$ small and $-\varepsilon\notin \sigma(A)$.
\end{list}

Consider the solutions of the following equation
\begin{equation}\label{s-d-1}
Az=\Phi'(z),\;z\in D(A).
\end{equation}
Let
\[
P_0=\int^\beta_{-\beta} dE_\lambda,\;P_+=\int^{+\infty}_\beta dE_\lambda,\;P_-=\int^{-\beta}_{-\infty} dE_\lambda,
\]
where $\{E_\lambda\}$ is the spectral resolution of $A$, and let
\[
\mathcal{H}_*=P_*\mathcal{H},\;*=0,\pm.
\]
Decompose the space $V$ as follows
\[
V=V_0\oplus V_-\oplus V_+,
\]
where $V_*=|A_\varepsilon|^{-1/2}\mathcal{H}_*,\;*=0,\pm$ and $A_\varepsilon=A+\varepsilon I$.

For each $u\in \mathcal{H}$, we have the decomposition
\[
u=u_++u_0+u_-,
\]
where $u_*\in \mathcal{H}_*,\;*=0,\pm$, let $z=z_++z_0+z_-$, with
\[
z_*=|A_\varepsilon|^{-1/2}u_*,\;*=0,\pm.
\]
Define a functional $f$ on $\mathcal{H}$ as follows:
\[
f(u)=\frac{1}{2}(\|u_+\|^2+Q_+\|u_0\|^2-Q_-\|u_0\|^2-\|u\|^2)-\Phi_\varepsilon(z),
\]
where $Q_+=\int^\infty_0dE_\lambda$, $Q_-=\int^0_{-\infty}dE_\lambda$, and $\Phi_\varepsilon(z)=\frac{\varepsilon}{2}\|z\|_\mathcal{H}+\Phi(z)$.

The Euler equation of this functional is the system
\begin{equation}\label{s-d-2}
u_\pm=\pm|A_\varepsilon|^{-1/2}P_\pm\Phi'_\varepsilon(z),
\end{equation}
\begin{equation}\label{s-d-3}
Q_\pm u_0=\pm|A_\varepsilon|^{-1/2}Q_\pm P_0\Phi'_\varepsilon(z).
\end{equation}
Thus $z=z_++z_0+z_-$ is a solution of \eqref{s-d-1} if and only if $u=u_++u_0+u_-$ is a critical point of $f$.
The implicit function can be applied, yielding a solution $z_\pm(z_0)$ for fixed $z_0\in V_0$, such that $z_\pm\in C^1(V_0,V_\pm)$.
Since dim$V_0$ is finite, all topologies on $V_0$ are equivalent, we choose  $\|\cdot\|_{\mathcal{H}}$ as it norm. We have
\[
u_\pm(z_0)=|A_\varepsilon|^{1/2}z_\pm(z_0)\in C^1(\mathcal{H}_0,\mathcal{H}),
\]
which solves the system \eqref{s-d-2}.

Let
\[
a(z_0)=f(u_+(z_0)+u_-(z_0)+u_0(z_0)),
\]
where $u_0(z_0)=|A_\varepsilon|^{1/2}z_0$ and let $z_0=x$, we have
\[
a(x)=\frac{1}{2}(A(z(x),z(x)))-\Phi(z(x)),
\]
where $z(x)=\xi(x)+x$, $\xi(x)=z_+(x)+z_-(x)\in D(A)$.
Then, we have the following theorem duo to Amann and Zehnder\cite{Amann}, Chang\cite{Chang-1993} and Long\cite{Long3}.
\begin{thm}\label{s-p-r}
Under the assumption (1), (2), (3), there is a one-one correspondence
\[
x\mapsto z=z(x)=z_+(x)+z_-(x)+x,
\]
between the critical points of the $C^2$-function $a\in C^2(\mathcal{H}_0,\mathbb{R})$ with the solutions of the operator equation
\[
Az=\Phi'(z),\;z\in D(A).
\]
Moreover, the functional $a$ satisfies
\[
a(x)=\frac{1}{2}(A(z(x),z(x)))-\Phi(z(x)),
\]
\[
a'(x)=A(z(x))-\Phi'(z(x))=Ax-P_0\Phi'(z(x)),
\]
\[
a''(x)=[A-\Phi''(z(x))]z'(x)=AP_0-P_0\Phi''(z(x))z'(x).
\]
\end{thm}
Since $\mathcal{H}_0$ is a finite dimensional space,  for every critical point $x$ of $a$ in $\mathcal{H}_0$, the Morse index and nullity  are finite, we denote them
by $(m^-_a(x), m^0_a(x))$.

Now, let the Hilbert space $\mathcal{H}$ be $L^2(\mathbb{R},\mathbb{R}^{2N})$, and the operator $A$ be $\tilde{F}=-J\frac{d}{dt}+L$,
$\Phi(u)=\int^\infty_{-\infty}R(t,u)$. Then we have $V=E$. For $R\in C^2(\mathbb{R}\times \mathbb{R}^{2N},\mathbb{R})$ and
$|\nabla^2_zR|\leq C_R,\;\forall (t,z)\in\mathbb{R}\times \mathbb{R}^{2N}$, let $-\alpha=\beta\geq2(C_R+1)$ and $\beta\not\in\sigma(\tilde{F})$,
we have $A$ and $\Phi$ satisfying the above conditions. Thus from Theorem \ref{s-p-r}, we can solve our
problems on the finite dimensional space.
 Similar to  Lemma 2.2 and Remark 2.3 in \cite{Liu-Su-Wang-adv}, we have the following  estimates.
\begin{lem}
Assume that $R\in C^2(\mathbb{R}\times \mathbb{R}^{2N},\mathbb{R})$, $|\nabla^2_zR|\leq C_R,\;\forall (t,z)\in\mathbb{R}\times \mathbb{R}^{2N}$
 and $\nabla_zR(t,0)\equiv0$, then we have
$$ \|u^\pm(x)\|_{L^2}\leq\frac{2\sqrt{\beta}(C_R+1)}{\beta-2C_R-3\varepsilon}\|x\|_{L^2},\forall x\in \mathcal{H}_0.$$
Moreover, we have
$$ \|(u^{\pm})'(x)\|_{L^2} \rightarrow 0,\beta \rightarrow \infty.$$
\end{lem}
\noindent{\it Proof.} Note that
$$u^\pm(x)=\pm |A_\varepsilon|^{-1/2}P_{\pm} \Phi'_\varepsilon(z^++z^-+x).$$
From $\nabla_zR(t,0)=0$, $|\nabla^2_zR|\leq C_R$, we have  $\Phi'(0)=0$ and $\|\Phi'(z)\|_{L^2}\leq C_R\|z\|_{L^2}$.  Since
$\||A_\varepsilon|^{-1/2}P_{\pm}\|\leq\frac{1}{\sqrt{\beta-\varepsilon}}$, we have
\begin{align}
\|u^\pm(x)\|_{L^2}&\leq\frac{1}{\sqrt{\beta-\varepsilon}}\|\Phi'(z^++z^-+x)+\varepsilon(z^++z^-+x)\|_{L^2}\nonumber\\
            &\leq\frac{C_R+\varepsilon}{\sqrt{\beta-\varepsilon}}\|z^++z^-+x\|_{L^2}\nonumber\\
            &\leq\frac{C_R+\varepsilon}{\sqrt{\beta-\varepsilon}}\left(\frac{\|u^+(x)\|_{L^2}}{\sqrt{\beta}}
            +\frac{\|u^-(x)\|_{L^2}}{\sqrt{\beta}}+\|x\|_{L^2}\right).
\end{align}
Therefore,
$$ \|u^+(x)\|_{L^2}+\|u^-(x)\|_{L^2}\leq\frac{2\sqrt{\beta}(C_R+\varepsilon)}{\beta-2C_R-3\varepsilon}\|x\|_{L^2}.$$
Next, since
$$ (u^\pm)'(x)=\pm |A_\varepsilon|^{-1/2}P_{\pm} \Phi'_\varepsilon(z^++z^-+x)((z^+)'(x)+(z^-)'(x)+I),$$
where $I$ is the identity map on $\mathcal{H}_0$, we have
$$\|(u^+)'(x)\|_{L^2}+\|(u^-)'(x)\|_{L^2}\leq\frac{2\sqrt{\beta}(C_R+\varepsilon)}{\beta-2C_R-3\varepsilon},\forall x\in \mathcal{H}_0.$$

\begin{rem}\label{Remark-z-estimate}
{\rm For $z(x)$, we also have that there is a constant $C>0$
dependent of $C_R$, but independent of $\beta$, such that
$$\|z^{\pm}(x)\|_{V}\leq\frac{C}{\sqrt{\beta}}\|x\|_{L^2},\|z'^{\pm}(x)\|_{V}\leq\frac{C}{\sqrt{\beta}},\forall x\in
\mathcal{H}_0.$$
 }
\end{rem}

If $R$ satisfies the condition $(R_1)$, then for any homoclinic orbit $z$ of {\rm (HS)}, $\nabla^2_zR(\cdot,z)\in\mathcal{B}$,
and hence we have  the associated index pair ($\mu_F(KB)$, $\upsilon_F(KB)$). For notation simplicity, in what follows, we set
\[
\mu_F(z)=\mu_F(K\nabla^2_z(R(t,z))),\]
 and
 \[
 \upsilon_F(z)=\upsilon_F(K\nabla^2_z(R(t,z))).
 \]

\begin{thm}\label{rm-sd}
Let  $R\in C^2(\mathbb{R}\times \mathbb{R}^{2N},\mathbb{R})$ satisfying
$|\nabla^2_zR|\leq C_R,\;\forall (t,z)\in\mathbb{R}\times \mathbb{R}^{2N}$ and $\nabla_zR(t,0)\equiv0$.
For each critical point $x$ of $a$ in $\mathcal{H}_0$, $z(x)$ is a  homoclinic orbit  of {\rm (HS)} and   we have
\begin{align}
m^-_a(x)&=\dim(E^-(\mathcal{H}_0))+\mu_F(K\nabla^2_z R(t,z(x)))=\dim(E^-(\mathcal{H}_0))+\mu_F(z(x)),\\
 m^0_a(x)&=\upsilon_F(K\nabla^2_z R(t,z(x))=\upsilon_F(z(x)),
\end{align}
where $\dim(E^-(\mathcal{H}_0))$ is the dimension of the space $\int^{0^-}_{-\beta}dE_\lambda(\mathcal{H}_0)$.
\end{thm}

This theorem shows the relations between  the relative Morse index and the Morse index of the saddle point reduction,
it will play an important role in the proof of our main results. The proof of this theorem will be postponed in the next subsection
where the notion of spectral flow will be used.
\subsection{The relationship between $\mu_F(T)$, spectral flow and the Morse index of saddle point reduction}
It is well known that the concept of spectral flow was first introduced by Atiyah, Patodi and Singer in \cite{Atiyah-Patodi-Singer-1976},
 and then extensively studied in \cite{Cappell-Lee-Miller-1994,Floer-1988,Robbin-Salamon-1993,Robbin-Salamon-1995,Zhu-Long-cam}.
 Here, we give a brief introduction of the spectral flow as introduced in \cite{Chen-Hu-2007}. Let $\mathcal{H}$ be a separable Hilbert space as defined before,
 and  $\{\mathcal{F}_\theta|\theta\in [0,1]\}$ be a continuous path of self-adjoint Fredholm operators on the Hilbert space $\mathcal{H}$.
 The spectral flow of $\mathcal{F}_\theta$ represents the net change in the number of negative eigenvalues of $\mathcal{F}_\theta$ as $\theta$ runs from 0 to 1, where the counting follows from the rule that each negative eigenvalue crossing to the positive axis contributes $+1$ and each positive eigenvalues crossing to the negative axis contributes $-1$, and for each crossing the multiplicity of eigenvalue is taken into account. In the calculation of spectral flow, a crossing operator introduced in \cite{Robbin-Salamon-1995} will be used. Take a $C^1$ path $\{\mathcal{F}_\theta|\theta
 \in[0,1]\}$ and let $\mathcal{P}_\theta$ be the projection from $\mathcal{H}$ to $\mathcal{H}^0(\mathcal{F}_\theta)$. When eigenvalue crossing occurs at $\mathcal{F}_\theta$, the operator
 \begin{equation}
 \mathcal{P}_\theta\frac{\partial}{\partial\theta}\mathcal{F}_\theta\mathcal{P}_\theta:\mathcal{H}^0(\mathcal{F}_\theta)\to \mathcal{H}^0(\mathcal{F}_\theta)
 \end{equation}
is called a crossing operator, denoted by $C_r[\mathcal{F}_\theta]$. As mentioned in \cite{Robbin-Salamon-1995}, an eigenvalue crossing at $\mathcal{F}_\theta$ is said to be regular if the null space of $C_r[\mathcal{F}_\theta]$ is trivial. In this case, we define
\begin{equation}
{\rm sign\;} C_r[\mathcal{F}_\theta]=\dim \mathcal{H}^+(C_r[\mathcal{F}_\theta])-\dim \mathcal{H}^-(C_r[\mathcal{F}_\theta]).
\end{equation}
A crossing occurs at $\mathcal{F}_\theta$ is called simple crossing if dim $\mathcal{H}^0(\mathcal{F}_\theta)=1$.

As indicated in \cite{Zhu-Long-cam}, the spectral flow $Sf(\mathcal{F}_\theta)$ will remain the same after a small disturbance of $\mathcal{F}_\theta$, that is, $Sf(\mathcal{F}_\theta)=Sf(\mathcal{F}_\theta+\varepsilon id)$ for $\varepsilon>0$ and small enough, where $id$ is the identity map on $\mathcal{H}$. Furthermore, we can choose suitable $\varepsilon$ such that all the eigenvalue crossings occurred in $\mathcal{F}_\theta, 0\leq\theta\leq1$ are regular\cite{Robbin-Salamon-1995}. Thus, without loss of generality, we may assume  all the crossings are regular. Let $\mathcal{D}$ be the set containing all the points in $[0,1]$ at which the crossing occurs. The set $\mathcal{D}$ contains only finitely many points. The spectral flow of $\mathcal{F}_\theta$ is
\begin{equation}
Sf(\mathcal{F}_\theta,0\leq\theta\leq 1)=\displaystyle\sum_{\theta\in\mathcal{D}^*}{\rm sign}C_r[\mathcal{F}_\theta]-\dim \mathcal{H}^-(C_r[\mathcal{F}_0])+\dim \mathcal{H}^+(C_r[\mathcal{F}_1]),
\end{equation}
where $\mathcal{D}^*=\mathcal{D}\cap(0,1)$. In what follows, the spectral flow of $\mathcal{F}_\theta$ will be simply denoted by $Sf(\mathcal{F}_\theta)$ when the starting and end points of the flow are clear from the contents. And $P_{\mathcal{F}_\theta}$ will be simply denoted by $P_\theta$.
\begin{prop}\label{index-sp}(See \cite[Proposition 3]{Chen-Hu-2007}.)  Suppose that, for each $\theta\in[0,1]$, $\mathcal{F}_\theta-\mathcal{F}_0$ is a compact operator on $\mathcal{H}$, then
$$ind(P_0|_{\mathcal{H}^-{(\mathcal{F}_1)}})=-Sf(\mathcal{F}_\theta).$$
\end{prop}
\noindent Thus, from Definition \ref{relative fredholm index},
$$\mu_{\mathcal{F}_0}(\mathcal{T})=-Sf(\mathcal{F}_\theta,\;0\leq\theta\leq 1),$$
where $\mathcal{F}_\theta=\mathcal{F}-\theta \mathcal{T}$, $\mathcal{T}$ is a compact operator. More over, if $\sigma(\mathcal{T})\subset [0,\infty)$ and $0\notin\sigma_P(\mathcal{T})$,  from the definition of Spectral flow, we have
\begin{align}
\mu_{\mathcal{F}_0}(\mathcal{T})&=-Sf(\mathcal{F}_\theta,\;0\leq\theta\leq 1)\nonumber\\
&=\displaystyle\sum_{\theta\in[0,1)}\upsilon_\mathcal{F}(\theta \mathcal{T})\nonumber\\
&=\displaystyle\sum_{\theta\in[0,1)}\dim \mathcal{H}^0(\mathcal{F}-\theta \mathcal{T}).
\end{align}
The proof of Theorem \ref{rm-sd} is the direct consequence of  the above Proposition \ref{index-sp} and Theorem3.2 in \cite{Zhu-Long-cam}, so we omit it here.
\begin{rem}\label{Remark-2.2}
The case of $R^-_\infty$ can be transformed into the case of $R^+_\infty$. More concretely, the $R^-_\infty$ case follows from the $R^+_\infty$ case by applying to the  function $\tilde{R}(t,z)=-R(-t,z)$.  If $z(t)$ is a homoclinic solution of $\tilde{F}z(t)=\nabla_z R(t,z(t))$, let $\tilde{z}(t)=z(-t)$, it's easy to check that $\tilde{z}(t)$ is a homoclinic solution of $\tilde{F}\tilde{z}(t)=\nabla_z\tilde{ R}(t,\tilde{z}(t))$, and this is a one-one correspondence between the two systems. By the definition of spectral flow and its catenation property\cite{Zhu-Long-cam}, we have $\mu_F(-B_\infty(-t))-\mu_F(-B_0(-t))=\mu_F(B_0(t))-\mu_F(B_\infty(t))$. Thus, we only consider the case of $R^+_\infty$ from now on.
\end{rem}

\section{Proof of our main results}
 {\it Proof of Lemma \ref{Lp-embedding}.} Recall the operator $\tilde{F}=-\mathcal {J}\frac{d}{dt}+L(t)$, with domain
$D(\tilde{F})=H^1(\mathbb{R},\mathbb{R}^{2N})$ if $L(t)$ is bounded and $D(\tilde{F})\subset H^1(\mathbb{R},\mathbb{R}^{2N})$ if $L(t)$ is unbounded. $D(\tilde{F})$ is a Hilbert space equipped with the norm
$
||z||_{\tilde{F}}=||(I+|\tilde{F}|)z||_{L^2},\;\forall z\in D(\tilde{F}).
$
Recall the  Hilbert space $E=D(|\tilde{F}|^{1/2})$, with the inner product and norm by
$$(u,v)_E=(|\tilde{F}|^{1/2}u,|\tilde{F}|^{1/2}v)_2+(u,v)_2,$$
$$\|u\|_E=(u,u)^{1/2}_E,$$
where $(\cdot,\cdot)_{L^2}$ denotes the usual inner product on $L^2(\mathbb{R},\mathbb{R}^{2N})$.
From ($L_1$), there is a matrix $L_0$ such that  $P(L(t)-L_0)\geq 0$ for all $t\in\mathbb{R}$.
We have $L(t)=\tilde{F}-\tilde{F}_0+L_0$, with $\tilde{F}_0=-\mathcal {J}\frac{d}{dt}+L_0$ and $D(\tilde{F}_0)=H^1$. Thus, for any $z\in E$
\begin{align}\label{eq-Lp-estimate-1}
|(L(t)z,P^Tz)_{L^2}|&\leq |(\tilde{F}z,P^Tz)_{L_2}|+|(\tilde{F}_0z,P^Tz)_{L_2}|+|(L_0z,P^Tz)_{L_2}|\nonumber\\
                 &\leq c\|z\|^2_E.
\end{align}
 Let $K\subset E$ be a bounded set.
We will show that $K$ is precompact in $L^p$ for $1\leq p\in (2/(1+\alpha), \infty)$. We derive the proof into three steps.

\noindent{\it Step 1. The case of  $p=2$.}  For $R>0$, from ($L_1$) and \eqref{eq-Lp-estimate-1} we have
\begin{align}\label{eq-Lp-estimate-2}
\int_{|t|>R}|z|^2&\leq c|R|^{\alpha-2}\int_{|t|>R}\langle L(t)z,P^Tz\rangle_{\mathbb{R}^{2N}}\nonumber\\
                 &\leq c|R|^{\alpha-2}\|z\|^2_E.
\end{align}
For any $\varepsilon>0$, from \eqref{eq-Lp-estimate-2}, we can choose $R_0$ large enough, such that
\beq\label{eq-Lp-estimate-3}
\int_{|t|>R_0}|z|^2<\frac{\varepsilon^2}{4}, \forall z\in K.
\eeq
On the other hand, by the definition of $\|\cdot\|_E$, we have
\beq\label{eq-Lp-estimate-4}
\int_{|t|\leq R_0}|z|^2\leq \|z\|_E\leq C, \forall z\in K.
\eeq
Thus, by the Sobolev compact embedding theorem there exist $z_1,z_2,\cdots,z_m\in K$, such that for any $z\in K$ there is $z_i$ satisfying
\beq\label{eq-Lp-estimate-5}
\|z-z_i\|^p_{L^p((-R_0,R_0),\mathbb{R}^{2N})}<\frac{\varepsilon^2}{2}.
\eeq
From \eqref{eq-Lp-estimate-3} and \eqref{eq-Lp-estimate-5}, we have $\|z-z_i\|_{L^2}<\varepsilon$, thus, $K$ has a finite $\varepsilon-$net in $L^2$, so the embedding $E\hookrightarrow L^2$ is compact.

\noindent{\it Step 2.  The case of  $p>2$.} Since $E$ is continuously embedded in $H^{1/2}$, hence by the Sobolev  embedding theorem, $E$ is continuously embedded in $L^p,\;\forall p>2$.
For any $p>2$, by the H\"{o}lder inequality we have
$$
\int_{\mathbb{R}}|z|^p\leq\|z\|_{L^2}\|z\|^{p-1}_{L^{2(p-1)}}\leq C\|z\|_{L^2}\|z\|^{p-1}_E,
$$
 thus, the embedding $E\hookrightarrow L^p$ is compact, $\forall p>2$.

\noindent{\it Step 3. The case of  $1\leq p\in (2/(1+\alpha),2)$.} First, we have $\frac{\alpha}{2-p}\cdot p>1$,  so we can choose $\alpha_p$ satisfying $\alpha_p\in (0,\alpha)$ and $\frac{\alpha_p}{2-p}\cdot p>1$. Denote by  $r=\frac{\alpha_p}{2-p}$. For $R>0$  and $z\in E$,
denote by $E^1_R(z)=\{t;|t|\geq R \;{\rm and}\; |t|^r|z(t)|>1\}$ and $E^2_R(z)=\{t;|t|\geq R \;{\rm and}\; |t|^r|z(t)|\leq1\}$.
Then, from \eqref{eq-Lp-estimate-1},
\begin{align}\label{eq-Lp-estimate-6}
\int_{E^1_R(z)}|z|^p&=\int_{E^1_R(z)}(|t|^r|z|)^p|t|^{-rp}\nonumber\\
                    &\leq\int_{E^1_R(z)}|z|^2|t|^{\alpha_p}\nonumber\\
                    &\leq \frac{c}{|R|^{\alpha-\alpha_p}}|(L(t)z,P^Tz)_{L^2}|\nonumber\\
                    &\leq \frac{c}{|R|^{\alpha-\alpha_p}}\|z\|^2_E,
\end{align}
and so
\begin{align}\label{eq-Lp-estimate-7}
\int_{|t|\geq R}|z|^p&=\int_{E^1_R(z)}|z|^p+\int_{E^2_R(z)}|z|^p\nonumber\\
                     &\leq \frac{c}{|R|^{\alpha-\alpha_p}}\|z\|^2_E+\frac{2}{(rp-1)R^{rp-1}},\;\forall z\in E.
\end{align}
Let $K\subset E$ be a bounded set. For any $\varepsilon>0$, from \eqref{eq-Lp-estimate-7}, choose $R_0>0$ large enough, such that
\beq\label{eq-Lp-estimate-8}
\int_{|t|\geq R_0}|z|^p<\frac{\varepsilon^p}{4}, \;\forall z\in K.
\eeq
On the other hand, by the Sobolev compact embedding theorem there are $z_1,z_2,\cdots, z_m\in K$, such that for any $z\in K$,
there exists $z_i$ satisfying
\beq\label{eq-Lp-estimate-9}
\|z-z_i\|^p_{L^p((-R_0,R_0),\mathbb{R}^{2N})}<\frac{\varepsilon^p}{2}
\eeq
From \eqref{eq-Lp-estimate-8} and \eqref{eq-Lp-estimate-9}, we have
$$
\|z-z_i\|_{L^p}< \varepsilon,
$$
that is to say $K$ has a finite $\varepsilon-$net in $L^p$, and the embedding $E\hookrightarrow L^p$ is compact. The proof of the lemma is compact.\endproof

Consider the homoclinic orbits of the linear Hamiltonian systems
\begin{equation}\label{linear systems}
\left\{\begin{array}{ll}\dot{z}(t)=JB(t)z(t), \forall t\in\mathbb{R},\\
       z(t)\to 0,\;|t|\to\infty.
       \end{array}
\right.
\end{equation}
where $z(t):\mathbb{R}\to \mathbb{R}^{2N}$,  $J=\left(\begin{matrix} 0 &-I_N\\I_N &0\end{matrix}\right)$ and $B(t)$ is a continuous symmetric matrix function.
Denote by $S$ the set of homoclinic orbits of linear systems \eqref{linear systems}, then $S$ is a linear subspace of $L^2(\mathbb{R},\mathbb{R}^{2N})$ and we have the following lemma.
\begin{lem}\label{astimate of nullity}
The dimension of the solution space  $S$ will be less than or equal to $N$. Thus for any homoclinic orbit $z(t)$ of (HS), if $R$ satisfies ($R_1$), we have
\[
0\leq\upsilon_F(z)\leq N.
\]
\end{lem}
\noindent{\it Proof.}  As usual, we define the symplectic groups on $\mathbb{R}^{2N}$ by
\[
Sp(2N)=\{M\in\mathcal{L}(\mathbb{R}^{2N}), |M^TJM=J\},
\]
where $\mathcal{L}(\mathbb{R}^{2N})$ is the set of all $2N\times 2N$ real matrices, $M^T$ denotes the transpose of $M$.
Let $W(t)$ be the fundamental solution of \eqref{linear systems}, then $W(t)$ is a path in $Sp(2N)$. Let $z(t)$ be a nontrivial homoclinic orbits of \eqref{linear systems}, that is to say $z(0)\neq 0$ and satisfies
\[
\left\{\begin{array}{ll}  z(t)=W(t)z(0),\\
\displaystyle\lim_{t\to\infty}W(t)z(0)=0.
 \end{array}\right.
\]
Denote by $S_0$ the subset of $\mathbb{R}^{2N}$ satisfying
\[
S_0=\{z\in\mathbb{R}^{2N}|\displaystyle\lim_{t\to\infty}W(t)z=0\},
\]
then we have $\dim{S}=\dim(S_0)$.
We claim that $Jz_0\not\in S_0$ if $z_0\in S_0$ and $z_0\neq 0$. We prove it indirectly, assume  $z_0, Jz_0\in S_0$ with $z_0\neq 0$, that is to say
\[
\displaystyle\lim_{t\to\infty}W(t)z_0=0,
\]
\[
\displaystyle\lim_{t\to\infty}W(t)Jz_0=0.
\]
Since $W(t)$ is a path in $Sp(2N)$, $W^T(t)JW(t)=J,\;\forall t\in\mathbb{R}$, thus
\begin{align}
0&=\displaystyle\lim_{t\to\infty}(JW(t)z_0, W(t)Jz_0)_{\mathbb{R}^{2N}}\nonumber\\
 &=-\displaystyle\lim_{t\to\infty}(z_0, W^T(t)J W(t)Jz_0)_{\mathbb{R}^{2N}}\nonumber\\
 &=(z_0, z_0)_{\mathbb{R}^{2N}},\nonumber
\end{align}
which contradicts $z_0\neq 0$. Since $J$ is an isomorphism on $\mathbb{R}^{2N}$, we have $\dim S_0\leq N$. And from the definition of $\upsilon_F(z)$ in the last part of subsection \ref{s-p-r}, we have complected the proof.\endproof

  Before the proof of Theorem \ref{main-result3}, we need the following lemma. Since $R$ satisfies condition ($R_1$), performing on $\rm (HS)$ the saddle point reduction. Choose a suitable number $\beta$, which is used in the projection for the saddle point reduction in
section 2.2. Let
\begin{equation}\label{}
P=\int^{\beta}_{-\beta}dE_\lambda,
\end{equation}
\begin{equation}\label{}
X=PL^2(\mathbb{R},\mathbb{R}^{2N}).
\end{equation}
By
Theorem \ref{s-p-r}, we have a functional
$a(x)$ with $x\in X$, whose critical points give rise to solutions of $\rm (HS)$.
\begin{lem}\label{lem-3.2}
(1) $a$ satisfies (PS) condition,\\
(2) $H_q(X,a;\mathbb{R})\cong\delta_{q,r}\mathbb{R},q=0,1,....$
for $-a\in \mathbb{R}$ large enough, where
$r=\dim(E^-(X))+\mu_F(KB_1).$
\end{lem}
\noindent{\it Proof.} Assume there is a sequence $\{x_n\}\subset X$, satisfying $a'(x_n)\to 0 (n\to \infty)$.
That is
\begin{equation}\label{eq-3.4}
\|Fz_n-K\nabla_zR(t,z_n)\displaystyle\|_E\to 0,
\end{equation}
 where $z_n=z(x_n)$ defined in section 2.1. Since $X$ is a finite dimensional space, and from the definition of $z_n$, it's enough to prove $\{z_n\}$ is bounded in $E$. For each $\varepsilon\in(0,1)$, define $C_n\in\mathcal{B}$ by
 \begin{equation}\label{C_n}
 C_n(t)=\left\{
 \begin{array}{ll}
\int^1_0\nabla^2_zR(t,sz_n)ds,\;|z_n(t)|\geq \frac{R_0}{\varepsilon},\\
 B_1(t),\;|z_n(t)|< \frac{R_0}{\varepsilon}.
 \end{array}
 \right.
 \end{equation}
It is easy to verify that $\{C_n\}$   satisfies
\[
B_1(t)-\varepsilon(B_1(t)+c\cdot I)\leq C_n(t)\leq B_2(t)+\varepsilon (c\cdot I-B_2(t)),\;\forall t\in\mathbb{R},
\]
where $c$ is the constant in condition $(R_1)$ and $I$ is the identity map on $\mathbb{R}^{2N}$.  Since $B_1\leq B_2$, $\mu_F(B_1)=\mu_F(B_2)$ and $\upsilon_F(KB_1)=\upsilon_F(KB_2)=0$, we can choose $\varepsilon$ small enough, such that for  each $n\in\mathbb{N}^+$, satisfying $\mu_F(KC_n)=\mu_F(KB_1)$ and $\upsilon_F(KC_n)=0$. Thus $F-KC_n$ is reversible on $E$ and there is a constant $\delta>0$, such that
\begin{equation}\label{eq-3.6}
\|(F-KC_n)z\|_E\geq\delta\|z\|_E,\;\forall z\in E, n\in\mathbb{N}^+.
\end{equation}
On the other hand, for $b\in(0,1)$, there is a constant $c>0$ depending on $b$, such that for each $n\in\mathbb{N}^+$,
\begin{equation}\label{eq-3.7}
|\nabla_zR(t,z_n(t))-C_nz_n(t)|\leq c|z_n(t)|^b, \forall t\in \mathbb{R}.
\end{equation}
Choose $b>\frac{1-\alpha}{1+\alpha}$ in \eqref{eq-3.7}, that is $1+b\in(\frac{2}{1+\alpha}, 2)$, we have
\begin{align}\label{eq-3.8}
\|(Fz_n-K\nabla_zR(t,z_n))-(F-KC_n)z_n\|_E^2&=\|K(\nabla_zR(t,z_n)-C_nz_n)\|_E^2\nonumber\\
                                            &\leq\|\nabla_zR(t,z_n)-C_nz_n\|_{L^2}^2\nonumber\\
                                            &\leq c\int_{\mathbb{R}}\frac{|\nabla_zR(t,z_n)-C_n z_n|}{|z_n|^{b}}|z_n|^{1+b} dt\nonumber\\
                                            &\leq c\|z_n\|_{L^{1+b}}^{1+b}.
\end{align}
As we claimed in the part of introduction,  in equations \eqref{eq-3.7} and \eqref{eq-3.8}
 the letter $c$ denotes different positive constants whose exact value is irrelevant.
Thus, from \eqref{eq-3.4}, \eqref{eq-3.6}, \eqref{eq-3.8} and Lemma \ref{Lp-embedding}, we have $\{z_n\}$ in bounded in $E$, and $a$ satisfies the ($PS$) conditions. And by Lemma 5.1 in
Chapter II of \cite{Chang-1993}, we have
$$H_q(X_,(a)_{\alpha};\mathbb{R})\cong\delta_{q,r}\mathbb{R},q=0,1,....,$$
 for $-\alpha\in \mathbb{R}$ large enough. $\hfill\Box$

 From Theorem\ref{rm-sd}, Lemma\ref{astimate of nullity} and Lemma\ref{lem-3.2},  Theorem\ref{main-result3} is a direct consequence of Theorem 5.1 and Corollary 5.2 in chapter II of \cite{Chang-1993}.

 In order to proof Theorem\ref{main-result1} and Theorem\ref{main-result2}, we need the following lemma
which is similar to  Lemma 3.4 in \cite{Liu-Su-Wang-calc} and Lemma 3.3 in  \cite{Liu-Su-Wang-adv}.

 \begin{lem}\label{Psim} Assume   $(R_1)$, $(R_0)$ and $(R^+_\infty)$ hold, then there exists a sequence of
 functions $R_k\in C^2(\mathbb{R}\times\mathbb{R}^{2N},\mathbb{R}),\;k\in\mathbb{N}$,
satisfying the following properties: \\
\rm (1)
 There exists an
increasing sequence of real numbers
$M_k\rightarrow\infty(k\rightarrow\infty)$ such that
\begin{equation}\label{eq2.1}
R_k(t,z)\equiv R(t,z),\;\,\forall\,  t\in\mathbb{R},\; |z|\leq M_k;
\end{equation}
\rm (2) For each $k\in\mathbb{N}$, there is a $C>0$ independent of $k$, such that
\begin{equation}\label{eq2.2'}
|\nabla^2_zR_k(t,z)|\leq C, \,\forall\,  t\in\mathbb{R},\; z\in\mathbb{R}^{2N},
\end{equation}
\begin{equation}\label{eq2.2}
\nabla^2_zR_k(t,z)\geq B_\infty,\,\forall\,  t\in\mathbb{R},\; |z|\geq R_0.
 \end{equation}
\rm (3) For each $k\in\mathbb{N}$, there exists some $C_k>0$ and a constant $\gamma$
with $\gamma I_{2N} >B_\infty$, $\nu_F(K\gamma I_{2N})=0$ such that
\begin{equation}\label{eq2.4}
|\nabla_zR_k(t,z)-\gamma z|<C_k, \;\,\forall\,  (t,z)\in \mathbb{R}\times\mathbb{R}^{2N},
\end{equation}
where $I_{2N}$ is the identity map on $\mathbb{R}^{2N}$.

\end{lem}
\noindent{\it Proof.}
Define $\eta:[0,\infty)\to \mathbb{R}$ by
$$\eta(s)=\left\{\begin{array}{ll} 0,&0\leq s< 1,\\
 \frac{2}{9}(s-1)^3-\frac{1}{9}(s-1)^4,&1\leq s<2,\\1-\frac{128}{9(12+s^2)},&2\leq s<\infty.
 \end{array}\right.$$
It's easy to see that $\eta\in C^2([0,\infty),\mathbb{R})$. Choose a sequence $\{M_k\}$ of positive
numbers such that $R_0<M_1<M_2<...<M_k<...\rightarrow\infty$ as
$k\rightarrow\infty$. For each $k\in \mathbb{N}$, let
$\eta_k(s)=\eta(\frac{s}{M_k})$ and
 \begin{equation}\label{eq2.8}
R_k(t,z)=(1-\eta_k(|z|))R(t,z)+\frac{\gamma}{2}\eta_k(|z|)|z|^2,\;k\in \mathbb{N}.
 \end{equation}
As in \cite{Liu-Su-Wang-calc,Liu-Su-Wang-adv}, we can  check that $R_k$ satisfies \eqref{eq2.1}--\eqref{eq2.4} for each $k\in \mathbb{N}$.
 \endproof

For each $k\in \mathbb{N}$, we consider the following problem
\[\left\{
\begin {array}{ll}
\tilde{F}z=\nabla_zR_k(t,z),\\
z(t)\to 0, z'(t)\to 0, t\to\infty
 \end{array}
 \right.\eqno(HS)_k\]
 where $R_k$ is given in Lemma \ref{Psim}.  Performing on $\rm (HS)_k$
the saddle point reduction. We choose the number $\beta$
which is used in the projection for the saddle point reduction in
section 2.2. First we choose
$$\beta>\text{max}\{2(C+1),2(\gamma +1)\}, {\rm and}\; \beta \not \in \sigma(A_0).$$
Let
\begin{equation}\label{eq-P}
P_{\beta}=\int^{\beta}_{-\beta}dE_\lambda,
\end{equation}
\begin{equation}\label{eq-X}
X_\beta=P_{\beta}L^2(\mathbb{R},\mathbb{R}^{2N}).
\end{equation}
 Thus for each $k$ and such a $\beta$ fixed, by
Theorem \ref{s-p-r}, we have a functional
$$a_{k,\beta}(x),\;x\in X_\beta,$$
whose critical points give rise to solutions of $\rm (HS)_k$.  Similarly
we have a functional
$$a_{\gamma,\beta}(x),\;x\in X_\beta,$$
whose critical points give
rise to solutions of the following systems $\rm (HS)_\gamma$
\[\left\{
\begin {array}{ll}\tilde{F}z=\gamma z,\\
z(t)\to 0, z'(t)\to 0, t\to\infty.
\end{array}
\right.\eqno{\rm (HS)_\gamma}\]
 For notational simplicity, we denote $a_k,\;a_\gamma$ for $a_{k,\beta}$ and $a_{\gamma, \beta}$. Define
 \[
 \Phi_k(z)=\int^\infty_{-\infty}R_k(t,z).
 \]
 For the functional $a_k$, similar to Lemma \ref{lem-3.2}, we have the following lemma.
\begin{lem}\label{ps-condition}
{\it\text{ }\\
  (1) $a_k$ satisfies (PS) condition,  the critical point set of $a_k$ is compact,\\
(2) $H_q(X_{\beta},(a_k)_{\alpha_k};\mathbb{R})\cong\delta_{q,r_\beta}\mathbb{R},q=0,1,....$
for $-\alpha_k\in \mathbb{R}$ large enough, where
$r_\beta=\dim(E^-(X_\beta))+\mu_F(K\gamma I_{2N}).$}
\end{lem}
 \noindent{\it Proof.} The proof is similar to Lemma \ref{lem-3.2}.
From Theorem \ref{rm-sd},  we have
\begin{align}\label{eq-ps1}
 \hspace{20mm}\|a'_k(x)-a''_\gamma(0)x\|_{L^2}
&=\|P_{\beta}K(\nabla_z\Phi_k(z_k(x))-\gamma x)\|_{L^2}\nonumber\\
&\leq\|P_{\beta}(\nabla_zR_k(z_k(x))-\gamma z_k(x))\|_{L^2}\nonumber\\
&\leq\|(\nabla_zR_k(z_k(x))-\gamma z_k(x))\|_{L^2}.
\end{align}
Similar to \eqref{eq-3.7}, we have for $b\in(0,1)$, there is some $c>0$, such that
\[
|\nabla_zR_k(t,z_k)-\gamma z_k|\leq c|z_k(t)|,\;\forall t\in\mathbb{R},
\]
Choose $b\in(\frac{1-\alpha}{1+\alpha},1)$, similar to \eqref{eq-3.8}, we have
\begin{align}\label{eq-ps2}
\|\nabla_zR_k(t,z_k)-\gamma z_k\|^2_{L^2}
&\leq c\|z_k\|^{1+b}_{L^{1+b}}.
\end{align}
From Lemma \ref{Lp-embedding}, Remark \ref{Remark-z-estimate} and  equation \eqref{eq-ps2},
\begin{align}\label{eq-ps3}
\|\nabla_zR_k(t,z_k)-\gamma z_k\|^2_{L^2}&\leq c\|z_k\|^{1+b}_E\nonumber\\
                                         &\leq c(\beta)\|x\|^{1+b}_{L^2}.
\end{align}
From \eqref{eq-ps1} and \eqref{eq-ps3},we have
\begin{equation}\label{eq-ps4}
\|a'_k(x)-a''_\gamma(0)x\|^2_{L^2}\leq c(\beta)\|x\|^{1+b}_{L^2}.
\end{equation}
 Now, for each $k\in\mathbb{N}$, we assume  $\{x_m\}\subset X_{\beta}$
satisfying $\|a'_k(x_m)\|\rightarrow 0$. By $\nu_F(K\gamma I_{2N})=0$, we have  $a''_\gamma(0)$ is invertible on $X_{\beta}$, since $b<1$
 the sequence  $\{x_m\}$ must be bounded. Thus the (PS) condition for $a_k$ holds. From the same reason,  we have the compactness of the critical point set of $a_k$.  And by Lemma 5.1 in
Chapter II of \cite{Chang-1993}, we have
$$H_q(X_{\beta},(a_k)_{\alpha_k};\mathbb{R})\cong\delta_{q,r_\beta}\mathbb{R},q=0,1,....,$$
 for $-\alpha_k\in \mathbb{R}$ large enough. $\hfill\Box$
\begin{lem}\label{uniform-boundary-of-solutions}
There exist $c>0$, such that for any $k\in \mathbb{N}$, and $z\in L^2$ satisfies the systems
$\rm (HS)_k$,
if $\mu_F(z)\leq \mu_F(KB_\infty)-1$,  we have
$\|z\|_{L^\infty}\leq c$.
\end{lem}
\noindent{\it Proof.}
 We prove it indirectly. Assume there exist $R_k$, $z_k$, satisfies the conditions, and $\|z_k\|_{L^\infty}\to\infty$,
that is $\|z\|_{\tilde{F}}\to\infty$.
Since $|\nabla_zR_k(t,z)|<c|z|,\;\forall t\in \mathbb{R},z\in\mathbb{R}^{2N}$, we have $\|z_k\|_{\tilde{F}}\leq c\|z_k\|_{L^2}$. Denote
$y_k=\frac{z_k}{\|z_k\|_{\tilde{F}}}$, then we have $y_k\to y$ in $L^2$ for some $y\in L^2$ with $\|y\|_{L^2}>0$, and
\[
\tilde{F}y_k=\frac{R'_k(t,z_k)}{\|z_k\|_{\tilde{F}}}.
\]
Then for any $r>0$, there exist $C_r>0$, satisfying
\begin{equation}\label{main-proof-3.32}
|\dot{y}_k(t)|\leq C_r|y_k(t)|,\;t\in I_r,
\end{equation}
where $I_r=[-r,r]$. Since $\|y\|_{L^2}>0$, there is a $r_0>0$, such that
\[\|y\|_{L^2(-r,r)}>\frac{1}{2}\|y\|_{L^2}>0,\;\forall r>r_0.
\]
Then from the similar argument in \cite{Liu-Su-Wang-adv}, there is a subsequence we may assume  $\{y_k\}$ converges in uniform norm to $y$,
and $y(t)\neq 0,\forall t\in I_r$. Therefor  $|z_k(t)|\to \infty$ uniformly on $I_r$, and there is $K(r)$ depending on $r$,
such that $|z_k(t)|\geq R_0$, for any $t\in I_r$ and $k\geq K(r)$.

Performing the saddle point reduction on the following systems
\[\left\{
\begin {array}{ll}\tilde{F}z=B_\infty z,\\
z(t)\to 0, z'(t)\to 0, t\to\infty.
\end{array}
\right.\eqno{\rm(HS)_\infty}\]
For $\beta$ large enough, we have the functional $a_{\infty,\beta}$ (denote by $a_\infty$ for simplicity) and the function $z(x)$, since $\upsilon_F(KB_\infty)=0$, we have the following decomposition
\[
X_\beta=X^+_\beta+X^-_\beta,
\]
where $a_\infty''(0)$ is positive definite on $X^+_\beta$ and negative definite on $X^-_\beta$.
From Remark\ref{Remark-z-estimate}, and $\upsilon_F(KB_\infty)=0$, there exists $\alpha>0$, such that for $\beta$ large enough
\begin{equation}
((\tilde{F}-B_\infty(t))x,x)_{L^2}\leq -\alpha\|x\|^2_{L^2},\;\forall x\in X^-_\beta.
\end{equation}
 From the uniform boundary of  $\nabla^2_zR_k(t,z)$ and Remark \ref{Remark-z-estimate},  we can choose $\beta$ large enough, such that
\begin{equation}
\|(\nabla^2_zR_k(t,z_k)z_k'^{\pm}(x_k)x,x\|_{L^2}\leq \frac{\alpha}{4}\|x\|_{L^2},\; \forall x\in L^2,
\end{equation}
where $x_k=P_\beta z_k$,  $z_k(x_k)=z_k(t)$ defined in Theorem\ref{s-p-r}.
 Choose $\varepsilon >0$ small enough and $\varepsilon<\frac{\alpha}{4}$, such that
$\mu_F(KB_\infty)=\mu_F(K(B_\infty-\varepsilon\cdot Id))$. Since $X^-_\beta$ is finite dimensional space, choose $r$ large enough, such that
\begin{equation}\label{main-proof-3.33}
((\nabla^2_zR_k(t,z_k)-B_\infty)x,x)_{L^2(I^c_r)}\geq-\varepsilon(x,x)_{L^2}, \;\forall x\in X^-_\beta,
\end{equation}
where $I^c_r=\mathbb{R}\setminus I_r$, and from the definition of $R_k$,
\begin{equation}\label{main-proof-3.34}
\nabla^2_zR_k(t,z_k(t))\geq B_\infty(t),\; t\in I_r,\,k\geq K(r),
\end{equation}
that is
\begin{equation}\label{main-proof-3.35}
((\nabla^2_zR_k(t,z_k)-B_\infty)x,x)_{L^2(I_r)}\geq 0, \;\forall x\in X^-_\beta,k\geq K(r).
\end{equation}
From \eqref{main-proof-3.33} and \eqref{main-proof-3.35},
\begin{equation}
(\nabla^2_zR_k(t,z_k)x,x)_{L^2}\geq ((B_\infty x,x)_{L^2}-\varepsilon(x,x)_{L^2},
\end{equation}
for $k$ large enough. Thus we have
\begin{align}
(a''_k(x_k)x,x)_{L^2}&=((\tilde{F}-\nabla^2_zR_k(t,z_k))x,x)_{L^2}-(\nabla^2_z(R_k(t,z_k))(z'^+_k+z'^-_k)x,x)_{L^2}\nonumber\\
                     &\leq((\tilde{F}-B_\infty)x,x)_{L^2}+\frac{\alpha}{2}\|x\|^2_{L^2}+\varepsilon\|x\|^2_{L^2}\nonumber\\
                     &\leq -\frac{\alpha}{4}\|x\|^2_{L^2}.
\end{align}
That is $m^-_{a_k}(x)\geq m^-_{a_\infty}(0)$, from Theorem\ref{rm-sd},
$m^-_{a_k}(x_k)=\dim(E^-(\mathcal{H}_0))+\mu_F(z_x(x_k))$, $m^-_{a_\infty}(0)=\dim(E^-(\mathcal{H}_0))+\mu_F(KB_\infty)$,
thus $\mu_F(z_x)\geq \mu_F(KB_\infty)$, which contradicts the assumption.\endproof

\noindent{\it Proof of Theorem \ref{main-result1}.}  As claimed in Remark \ref{Remark-2.2}, we can only consider the case of $(R^+_\infty)$.
 Note that $z=0$ is a critical point of $a_k$, the morse index of $0$ for $a_k$ is $m^-_{a_k}(0)=\dim(E^-(\mathcal{H}_0))+\mu_F(KB_0)$,
since $\gamma\cdot I_{2N}>B_\infty$, we have
\begin{equation}
\mu_F(K\gamma\cdot I_{2N})\geq \mu_F(KB_\infty)\geq \mu_F(KB_0).
\end{equation}
From  proposition (2) in Lemma\ref{ps-condition}, use  the $(m^-_{a_k}(0))^{th}$ and $(m^-_{a_k}(0)+1)^{th}$ Morse inequalities,
 $a_k$ has a nontrivial critical point $x_k$ with it morse index $m^-_{a_k}(x_k)\leq m^-_{a_k}(0)+1$,
 that is $\mu_F(z_k)\leq \mu_F(KB_0)+1\leq \mu_F(KB_\infty)-1$, then from Lemma\ref{uniform-boundary-of-solutions},
we have $\{z_k\}$ is bounded in $L^\infty$. Thus $z_k$ is a nontrivial solution of (HS) for $k$ large enough.\endproof

The proof of Theorem \ref{main-result2} is similar to the proof of Theorem
\ref{main-result1}.  Instead of Morse theory we make use of minimax arguments for multiplicity of critical
points.

Let $X$ be a Hilbert space and assume $\phi\in C^2(X,\mathbb{R})$ is
an even functional, satisfying the (PS) condition and $\phi(0)=0$.
Denote $S_c=\{u\in X|\; \|u\|=c\}$.

 \begin{lem}\label{lem2.6}（See \cite[Corollary 10.19]{G}.) Assume $Y$ and $Z$ are
 subspaces of $X$ satisfying $\dim Y=j>k={\rm codim} Z$. If there exist
 $R>r>0$ and $\alpha>0$ such that
 $$\inf \phi (S_r\cap Z)\geq \alpha, \;\sup \phi (S_R\cap Y)\leq 0,$$
 then $\phi$ has $j-k$ pairs of nontrivial critical points $\{\pm x_1,\pm x_2,...,\pm
 x_{j-k}\}$,
so that $\mu (x_i)\leq k+i$, for $i=1,2,...j-k$.
\end{lem}
 First, we consider the case of $(R^+_\infty)$, since $R$ is even, we have $R_k$ is also even, and satisfies
Lemma \ref{Psim}. Let $Y=X^-_\beta$, and
$Z$ the positive space of $a''_k(0)$ in $X_\beta$, and we have
dim $Y=E^-(X_\beta)+\mu_F(KB_\infty)$, codim$ Z=E^-(X_\beta)+\mu_F(KB_0)$, dim$ Y>$codim$Z$.  So $a_k$ has $l:=\mu_F(KB_\infty)-\mu_F(KB_0)$ pairs of
 nontrivial critical points $$\{\pm x_1,\pm x_2,...,\pm x_{l}\},$$ and
 $l-1$ pairs of them satisfy
 \begin{equation}\label{eq2.48}
m^-(x_i)\leq \mu_F(KB_0)+i<\mu_F(KB_\infty),\;\;i=1,2,...,l-1.
 \end{equation}
Then we can complete the proof. In order to prove the case of $(R^-_\infty)$, we need the following lemma.
 \begin{lem}(See \cite[Corollary II 4.1]{Chang-1993}.)
 Assume $Y$ and $Z$ are subspaces of $X$ satisfying ${\rm dim} Y=j>k={\rm codim} Z$. If there exist $r>0$, and $\alpha>0$
 such that
 $$\inf \phi(Z)>-\infty, \sup \phi(S_r\cap Y)\leq -\alpha,$$
 then $\phi$ has $j-k$ pairs of nontrivial critical points $\pm u_1, \pm u_2, \cdots, \pm u_{j-k}$ so that
 $\mu(u_i)+\nu(u_i)\geq k+i-1$ for $i=1,2, \cdots, j-k$.
 \end{lem}
The proof is similar to the case of $(R^+_\infty)$, we omit it here.

\end{document}